\documentclass{amsproc}
\usepackage{amsfonts,amsmath,amssymb,amscd,amsthm,mathrsfs}

\newtheorem{theorem}{Theorem}[section]

\theoremstyle{definition}

\theoremstyle{remark}
\newtheorem{remark}[theorem]{Remark}

\numberwithin{equation}{section}

\begin{document}

\title[New Developments in Mean Curvature Flow]{New Developments in Mean Curvature Flow of Arbitrary Codimension Inspired By Yau Rigidity Theory}
\author{Li Lei, Hong-Wei Xu}



\thanks{Research supported by the National Natural Science Foundation of China, Grant Nos. 11531012, 11371315.}
\begin{abstract}
In this survey, we will focus on the mean curvature flow theory with
sphere theorems, and discuss the recent developments on the
convergence theorems for the mean curvature flow of arbitrary
codimension inspired by the Yau rigidity theory of submanifolds.
Several new differentiable sphere theorems for submanifolds are
obtained as consequences of the convergence theorems for the mean
curvature flow. It should be emphasized that Theorem 4.1 is an
optimal convergence theorem for the mean curvature flow of arbitrary
codimension, which implies the first optimal differentiable sphere
theorem for submanifolds with positive Ricci curvature. Finally, we
present a list of unsolved problems in this area.

\end{abstract}
\maketitle

\section{Huisken's classical theorems for MCF of hypersurfaces}
Let $(M,g)$ be a closed Riemannian $n$-manifold, and let $F_t
:M^n\rightarrow N^{n+q}$ be a one-parameter family of smooth
submanifolds immersed in an $(n+q)$-dimensional Riemannian manifold
$(N,h)$.
    We say that $M_t = F_t(M)$ is a solution of the mean curvature flow
    if $F_t$ satisfies
    \begin{eqnarray}
    \label{MCF}\left\{
    \begin{array}{ll}
    \frac{\partial}{\partial t}F(x,t)=H(x,t),\\
    F(x,0)=F_0(x),
    \end{array}\right.
    \end{eqnarray}
    where $F(x,t)=F_t(x)$, $H(x,t)$ is the mean curvature vector, and
    $F_{0}(M)$ is the initial submanifold immersed in $N$.

In 1984, Huisken first studied the mean curvature flow for compact
hypersurfaces in the Euclidean and proved the following convergence
theorem \cite{MR772132}.

\medskip

\noindent\textbf{Theorem 1.1.}
\emph{Let $M_0$ be an $n$-dimensional ($n \ge 2$) closed hypersurface in $ \mathbb{R}^{n+1}$. Assume that $M_0$ is uniformly convex. Then the mean curvature flow with initial value $M_0$
    has a smooth solution on a finite time interval $0 \le  t < T$,
    and converge to a single round point as $t \rightarrow T$.
}\medskip

Afterwards, Huisken \cite{MR837523} generalized Theorem 1.1 to the
case of the mean curvature flow of compact hypersurfaces in a
Riemannian manifold $N^{n+1}$ with bounded geometry, which is stated
as follows.
\medskip

\noindent\textbf{Theorem 1.2.}
\emph{Let $n \ge 2$ and $N^{n+1}$ be a complete Riemannian manifold
    which satisfies uniform bounds
    \begin{eqnarray*}
    -K_1\leq  \bar{K}_N \leq K_2,\,\,\,\,
    |\bar{\nabla}\bar{R}|\leq L^2,\ \  \,\,
    {\rm inj}(N)\geq i_N,\ \ \
    \end{eqnarray*}
    for nonnegative constants $K_1$, $K_2$, $L$ and positive constant
    $i_N$. Let $M_0$ be a closed hypersurface
    immersed in $N$, and suppose that $M_0$ satisfies the
    pinching condition
    \[|H| h_{ij} > n K_1 g_{ij} + \frac{n^2}{|H|}
     L g_{ij} .\]
    Then the mean curvature flow with initial value $M_0$
    has a smooth solution on a finite time interval $0 \le  t < T$,
    and converge to a single round point as $t \rightarrow T$.
}\medskip

Motivated by the rigidity theorem for hypersurfaces with constant
mean curvature in a sphere due to Okumura \cite{Okumura1}, Huisken
\cite{Huisken2} proved the following convergence theorem for the
mean curvature flow in a sphere.
\medskip

\noindent\textbf{Theorem 1.3.}
\emph{Let $n \ge 2$ and $\mathbb{F}^{n+1}(c)$  be a spherical space form with positive constant
    curvature $c$. Let $M_0$ be a closed hypersurface immersed in $\mathbb{F}^{n+1}(c)$, and suppose that $M_0$ satisfies the
    pinching condition
    \begin{eqnarray}
    |A|^2<\begin{cases}
    \frac{3}{4}|H|^2+\frac{4}{3}c, \ &n = 2, \\
    \frac{1}{n-1}|H|^2+2c, \ &n \geq 3.
    \end{cases}
    \end{eqnarray}
    Then one of the following holds:\\
    $(i)$  The mean curvature flow with initial value $M_0$ has a smooth solution $M_t$ on a finite time interval $0 \le  t < T$
    and the $M_t$'s converge uniformly to a single round point as $t \rightarrow T$.\\
    $(ii)$  The mean curvature flow with initial value $M_0$ has a smooth solution $M_t$ for all $0 \le  t < \infty$ and the $M_t$'s
    converge in the $C^\infty$-topology to a smooth totally geodesic hypersurface $M_\infty$.
}\medskip

During the past three decades, there are many progresses on the
theory of mean curvature flows. Most of the results focus on the
mean curvature flow of hypersurfaces. On the other hand, fruitful
results on the mean curvature flow of submanifolds of higher
codimension were obtained by several geometers. Mean curvature flow
of surfaces in 4-manifolds, Lagrangian mean curvature flow of higher
codimension, graphic mean curvature flow, mean curvature flow with
convex Gauss image were investigated by Chen, Li, Smoczyk, Wang, Xin
and others
\cite{ChenLi1,Smoczyk4,Wang1,Wang2,Wang3,Wang4,Wang5,Xin}.

In the present article, we will introduce the recent developments on
the convergence theorems for the mean curvature flow of arbitrary
codimension inspired by the Yau rigidity theory of submanifolds with
parallel mean curvature. Several new differentiable sphere theorems
for submanifolds are obtained as consequences of the convergence
theorems for the mean curvature flow. Finally, we present a list of
conjectures in this area.

\section{Yau rigidity theory of submanifolds and its developments}
More than forty years ago, Okumura \cite{Okumura,Okumura1}
investigated compact submanifolds with parallel mean curvature and
flat normal bundle in $\mathbb{S}^{n+q}$ whose squared norm of the
second fundamental form satisfies $|A|^2< 2+\frac{|H|^2}{n-1}$.
Meanwhile, Chen-Okumura \cite{ChenOku} proved a rigidity theorem for
compact submanifolds with parallel mean curvature in Euclidean
spaces under pinching condition $|A|^2< \frac{|H|^2}{n-1}.$ In 1975,
Yau \cite{Yau} established the rigidity theory of submanifolds with
parallel mean curvature, which includes the following important
results.
\medskip

\noindent\textbf{Theorem 2.1.} \emph{Let $M$ be an $n$-dimensional
compact submanifold with parallel mean
    curvature in $\mathbb{S}^{n+q}$ with $q>1$. If
    $|A|^2< n[3+n^{\frac{1}{2}}-(q-1)^{-1}]^{-1},$ then $M$ lies in
    a totally geodesic sphere $\mathbb{S}^{n+1}$.
}\medskip

\noindent\textbf{Theorem 2.2.} \emph{Let $M$ be an $n$-dimensional
oriented compact minimal submanifold
    in $\mathbb{S}^{n+q}$. If $K_{M}\geq\frac{q-1}{2q-1},$ then either
    $M$ is the totally geodesic sphere, the standard immersion of the
    product of two spheres, or the Veronese surface in $\mathbb{S}^4$.
}\medskip

Notice that Yau's pinching constant $\frac{q-1}{2q-1}$ is better
than the pinching constant $\frac{1}{2}$ given by Simon
\cite{Simon}. Moreover, Yau's pinching constant above is the best
possible in the case where $q=1$, or $n=2$ and $q=2$. Inspired by
the Yau rigidity theorem, Ejiri \cite{Ejiri} obtained a rigidity
theorem for compact and simply connected minimal submanifolds in a
sphere with $Ric_{M}\ge n-2.$  The Yau parameter method was first
introduced in the proof of the Yau rigidity theorem, which has many
applications in the study of the Chern conjecture and other rigidity
problems.

Many other important results were obtained by Yau \cite{Yau}. For
example, he obtained the geometric structure theorem for
non-negatively curved submanifolds with parallel mean curvature in
spaces forms, the classification theorem for surfaces with parallel
mean curvature in spaces forms, and the codimension reduction
theorem for submanifolds in a conformally flat manifold.

The Yau rigidity theory of submanifolds plays a very important role
in the study of geometry, topology and curvature flows of
submanifolds. During the past four decades, the Yau rigidity theory
has been developed by several geometers \cite{GXXZ}. In 1984, Cheng,
Li, and Yau \cite{Cheng} proved the following volume gap theorem for
minimal submanifolds.
\medskip

\noindent\textbf{Theorem 2.3.} \emph{Let $M$ be an $n$-dimensional
compact
    minimal submanifold  in $\mathbb{S}^{n+q}$. Then there exists an explicit positive
    constant $\varepsilon_n$ depending only on $n$ such that if
    $Vol(M)<Vol(\mathbb{S}^{n})+\varepsilon_n$, then $M$ is congruent to the great sphere $\mathbb{S}^{n}$.
}\medskip

The possible generalization of the Cheng-Li-Yau theorem will give a
proof of the higher dimensional version of the Willmore conjecture
on the total mean curvature for closed submanifolds in Euclidean
spaces.

In 1990, Xu \cite{Xu90} proved the generalized
Simons-Lawson-Chern-do Carmo-Kobayashi theorem for compact
submanifolds with parallel mean curvature in a sphere.
\medskip

\noindent\textbf{Theorem 2.4.} \emph{Let $M$ be an $n$-dimensional
oriented compact submanifold with
    parallel mean curvature in
    $\mathbb{S}^{n+q}$. If $|A|^2\leq C_0(n,q,|H|),$ then $M$ is either a
    totally umbilic sphere $\mathbb{S}^{n} ( n/ \sqrt{| H |^{2} +n^{2}}
    )$, a Clifford hypersurface in an $(n+1)$-sphere, or  the Veronese
    surface in $\mathbb{S}^{4} ( 2/ \sqrt{| H |^{2} +4} )$.}
\emph{Here the constant $C_0(n,q,|H|)$ is defined by
$$C_0(n,q,|H|):=\left\{\begin{array}{llll} \alpha(n,|H|),&q=1, \mbox{\ or\ }
q=2 \mbox{\ and\ } |H|\neq0,  &\\
\frac{n}{2-\frac{1}{q}}, & q\geq2  \mbox{\ and\ } |H|=0,  &\\
\min\Big\{\alpha(n,|H|),\frac{n^2+|H|^2}{(2-\frac{1}{q-1})n}+\frac{|H|^2}{n}\Big\},& q\geq 3 \mbox{\ and\ } |H|\neq0,
\end{array} \right.$$
\[\alpha (n,|H| ) =n+ \frac{n}{2 ( n-1 )} |H|^{2}- \frac{n-2}{2 ( n-1 )} \sqrt{|H|^{4}
    +4 ( n-1 ) |H|^{2}} .\]
}\medskip

In \cite{Li0,Li}, Li and Li improved Simons' pinching constant for
$n$-dimensional compact minimal submanifolds in $\mathbb{S}^{n+q}$
to $\max\{\frac{n}{2-1/q},\frac{2}{3}n\}$. Using Li-Li's matrix
inequality, Xu \cite{Xu} improved the pinching constant
$C_0(n,q,|H|)$ in Theorem 2.4 to
$$C(n,q,|H|)=\left\{\begin{array}{llll} \alpha(n,|H|),& q=1, \mbox{\ or\ } q=2 \mbox{\ and\ }
|H|\neq0,\\
\min\Big\{\alpha(n,|H|),\frac{2n}{3}+\frac{5}{3n}|H|^2\Big\},&\mbox{otherwise.\
}
\end{array} \right.$$
This is the best pinching constant for $n$-dimensional compact
submanifolds with parallel mean curvature in $\mathbb{S}^{n+q}$ up
to date. In fact, Theorem 2.4 is stronger than the rigidity results
in \cite{Alencar,Araujo,MR1289187}. In \cite{xuhwrigidity}, Xu
obtained an optimal rigidity theorem for complete submanifolds with
parallel mean curvature in hyperbolic spaces.

Let $\mathbb{F}^{n+q}(c)$ be an $(n+q)$-dimensional simply connected
space form with constant curvature $c$. Put
$$
\alpha (n,|H|,c ) =n c+ \frac{n}{2 ( n-1 )} |H|^{2}- \frac{n-2}{2 ( n-1 )} \sqrt{|H|^{4}
    +4 ( n-1 ) c |H|^{2}} .
$$
If $c>0,$ we have $\min_{|H|}\alpha(n,|H|,c)=2\sqrt{n-1} c.$

Using the techniques in stable currents and algebraic topology,
Shiohama and Xu \cite{Shiohama} proved the following optimal
topological sphere theorem.
\medskip

\noindent\textbf{Theorem 2.5.} \emph{Let $M^{n}(n\ge4)$ be an
oriented complete submanifold in $\mathbb{F}^{n+q}(c)$ with $c\geq
0$. Suppose that $\sup_{M}(|A|^2-\alpha(n,|H|,c))<0.$ Then $M$ is
homeomorphic to a sphere. }\medskip

After the work by Shen \cite{Shen} and Lin-Xia \cite{Lin}, Xu
\cite{Xu90} proved the following theorem.
\medskip

\noindent\textbf{Theorem 2.6.} \emph{Let $M$ be an $n$-dimensional
compact submanifold with parallel mean curvature in the space form
$\mathbb{F}^{n+q}(c) $ with  $c\geq 0$. Denote by $\mathring{A}$ the
traceless second fundamental form of $M$. If
$$\int_{M}|\mathring{A}|^n dM<C_1(n),$$ where $C_1(n)$ is an explicit
positive constant depending only on $n$, then $ M$ is a totally
umbilic sphere. }\medskip

Applying the Morse theory of submanifolds, Shiohama and Xu
\cite{Shiohama0} proved the following topological sphere theorem for
compact submanifolds.
\medskip

\noindent\textbf{Theorem 2.7.} \emph{Let $M$ be an $n$-dimensional
compact submanifold in the space form $\mathbb{F}^{n+q}(c) $ with
$c\geq0$. There exists a positive constant $C_2(n)$ depending only
on $n$ such that if
$$\int_{M}|\mathring{A}|^ndM< C_2(n),$$ then $M $ is homeomorphic to $\mathbb{S}^n$.
}\medskip

After the work due to Ejiri \cite{Ejiri}, Shen \cite{ShYB}, Li
\cite{Li2} and Xu-Tian \cite{XT}, Xu-Gu \cite{XG5} proved the
following pinching theorem for submanifolds with pinched Ricci
curvatures.
\medskip

\noindent\textbf{Theorem 2.8.} \emph{Let $M$ be an
$n(\geq3)$-dimensional oriented compact submanifold with  parallel
mean curvature in the space form $\mathbb{F}^{n+q}(c)$. If
$$Ric_{M}\geq(n-2)(c+\frac{|H|^2}{n^2}),$$ where  $n^2 c+|H|^2>0$, then $M$
is either a totally umbilic sphere, a Clifford hypersurface
$\mathbb{S}^{m}(r)\times \mathbb{S}^{m}(r)$ in an $(n+1)$-sphere
with $n=2m$, or $\mathbb{C}P^{2}(\frac{4}{3}c+\frac{1}{12}|H|^2)$ in
$\mathbb{S}^7(\frac{4}{\sqrt{16c+|H|^2}})$.  }\medskip

\begin{remark}
    When $n$ is even, the pinching condition in Theorem
    2.8 is optimal. Furthermore, Gu-Tian-Xu \cite{XG5,XT} obtained the sharp topological and differentiable sphere theorems for
    submanifolds with pinched Ricci curvatures.
\end{remark}

Using the DDVV inequality verified by Ge-Tang \cite{Ge} and Lu
\cite{Lu4} and the Yau parameter method, Gu-Xu \cite{GX} proved the
following rigidity theorem for minimal submanifolds in spheres.
\medskip

\noindent\textbf{Theorem 2.9.} \emph{Let $M$ be an $n$-dimensional
oriented compact minimal submanifold in the unit sphere
$\mathbb{S}^{n+q}$. If  $$K_{M}\geq\frac{q\cdot sgn(q-1)}{2(q+1)},$$
then $M$ is either a totally geodesic sphere, the standard immersion
of the product of two spheres, or the Veronese surface in
$\mathbb{S}^4$. Here $sgn(\cdot)$ is the standard sign function.
}\medskip
\begin{remark}
    When $2<q<n$, the pinching constant in Theorem 2.9 is better than
    ones given by Yau \cite{Yau} and Itoh \cite{Itoh2}.
\end{remark}

Combing Theorems 2.2, 2.9, and rigidity results in
\cite{GX,Itoh2,Shen2}, we present a general version of the Yau
rigidity theorem for submanifolds with parallel mean curvature in
spaces forms.
\medskip

\noindent\textbf{Generalized Yau rigidity theorem.} \emph{Let $M$ be
an $n$-dimensional oriented compact submanifold with parallel mean
curvature in $\mathbb{F}^{n+q}(c)$, where $n^2c+|H|^2>0$. Set
    $k(m,n)=\min\{m\cdot sgn(m-1), n\}.$ Then we have\\\\
$(i)$ if $|H|=0$ and
$$K_{M}\geq\frac{k(q,n)c}{2[k(q,n)+1]},$$ then M is either a
totally geodesic sphere, the standard immersion
of the product of two spheres, or the Veronese submanifold in $\mathbb{F}^{n+d}(c)$, where $d=\frac{1}{2}n(n+1)-1$;\\\\
$(ii)$ if $|H| \neq 0$ and $$K_M \ge
\frac{k(q-1,n)(c+\frac{1}{n^2}|H|^2)}{2[k(q-1,n)+1]},$$ then $M$ is
either a totally umbilical sphere, a Clifford hypersurface in an
$(n+1)$-sphere, a product of three spheres in an $(n+2)$-sphere, or
the Veronese submanifold in
$\mathbb{F}^{n+d}(c+\frac{1}{n^2}|H|^2)$, where
$d=\frac{1}{2}n(n+1)-1$. }
\medskip
\begin{remark}
Motivated by the generalized Yau rigidity theorem, Gu and Xu
\cite{XG6} verified a sharp differentiable sphere theorem for
submanifolds with positive sectional curvature.
\end{remark}

\section{MCF meets Ricci flow: positive sectional curvature}

Using the Ricci flow and stable currents
\cite{Brendle1,Brendle3,Hamilton,Lawson2}, Xu and Zhao \cite{XZ}
initiated the study of differentiable pinching problems for
submanifolds, and proved the following differentiable sphere
theorem.
\medskip

\noindent\textbf{Theorem 3.1.}
\emph{Let $M$ be an
$n(\geq4)$-dimensional oriented complete submanifold in the unit
sphere $\mathbb{S}^{n+q}$.
Then\\
$(i)$ if $n=4,5,6$ and $\sup_M(|A|^2-\alpha(n,|H|))<0$, then $M$ is diffeomorphic to $\mathbb{S}^{n}$;\\
$(ii)$ if $n\geq7$ and $|A|^2<2\sqrt{2}$, then $M$ is
diffeomorphic to  $\mathbb{S}^{n}$.}\medskip

In \cite{XG1}, Xu and Gu proved an optimal differentiable sphere
theorem for submanifolds in the space form $\mathbb{F}^{n+q}(c)$
with constant curvature $c\ge0$.
\medskip

\noindent\textbf{Theorem 3.2.} \emph{Let $M$ be an $n$-dimensional
oriented complete submanifold in $\mathbb{F}^{n+q}(c)$ with $c\ge0$.
If $\sup_M\Big(|A|^2- \frac{|H|^{2}}{n-1}-2c\Big)<0,$
 then $M$ is diffeomorphic to $\mathbb{S}^n$.}
\medskip

Meanwhile, Andrews-Baker \cite{Andrews-Baker} proved the following
convergence theorem for the mean curvature flow of arbitrary
codimension in Euclidean spaces, which implies a differentiable
sphere theorem for submanifolds in Euclidean spaces.
\medskip

\noindent\textbf{Theorem 3.3.} \emph{Let $M_0$ be an $n$-dimensional
compact submanifold immersed in $\mathbb{R}^{n+q}$. If $M_0$ has
$|H|\neq0$ everywhere and satisfies
\begin{eqnarray}
\label{pinch-cond00} |A|^2\leq\begin{cases}
\frac{4}{3n}|H|^2, \ &n = 2, 3, \\
\frac{1}{n-1}|H|^2, \ &n \geq 4,
\end{cases}
\end{eqnarray}
then the mean curvature flow has a smooth solution
$F:M\times[0,T)\rightarrow \mathbb{R}^{n+q}$ on a finite maximal
time interval, and $M_t$ converges uniformly to a point $q\in
\mathbb{R}^{n+q}$ as $t\rightarrow T$. The rescaled maps
$\widetilde{F}_t=\frac{F_t-x}{\sqrt{2n(T-t)}}$ converge in
$C^{\infty}$ to a limiting embedding $\widetilde{F}_T$ such that
$\widetilde{F}_T(M)$ is the unit $n$-sphere in some
$(n+1)$-dimensional subspace of $\mathbb{R}^{n+q}$.}\medskip

In \cite{GX1}, Gu and Xu obtained a refined version of Theorem 3.2.
Later, Baker \cite{Baker} proved a sharp convergence theorem for the
mean curvature flow of submanifolds in spheres.
\medskip

\noindent\textbf{Theorem 3.4.} \emph{Let $M_0$ be an $n$-dimensional
compact submanifold immersed in $\mathbb{S}^{n+q}$. If $M_0$
satisfies
\begin{eqnarray}
|A|^2\leq\begin{cases}
\frac{4}{3n}|H|^2+\frac{2(n-1)}{3}, \ &n = 2, 3, \\
\frac{1}{n-1}|H|^2+2, \ &n \geq 4,
\end{cases}
\end{eqnarray}
then either $M_t$ shrinks uniformly to a round point $p\in
\mathbb{S}^{n+q}$ as $t\rightarrow T<\infty$, or $M_t$ converges to
a totally geodesic sphere in $\mathbb{S}^{n+q}$ as $t$ tends to
infinity.}\medskip

In \cite{MR3078951}, Liu-Xu-Ye-Zhao  proved the following sharp
convergence theorem for the mean curvature flow of submanifolds in
hyperbolic spaces.
\medskip

\noindent\textbf{Theorem 3.5.} \emph{Let $F_0:M^{n}\rightarrow
\mathbb{F}^{n+q}(c)$ be an $n$-dimensional closed submanifold in a
hyperbolic space with constant curvature $c<0$. Assume $F_0$
satisfies
\begin{eqnarray}
\label{pinch-cond-H^N} |A|^2\leq\begin{cases}
\frac{4}{3n}|H|^2+\frac{n}{2}c, \ &n = 2, 3, \\
\frac{1}{n-1}|H|^2+2c, \ &n \geq 4.
\end{cases}
\end{eqnarray}
Then the mean curvature flow with $F_0$ as initial value converges
to a round point in finite time.
In particular, $M$ is diffeomorphic to
$\mathbb{S}^n$.}\medskip

More generally,  Xu and Gu \cite{XG1} proved the following
differentiable sphere theorem, which generalized the Brendle-Schoen
\cite{Brendle3} differentiable sphere theorem for manifolds with
strictly $1/4$-pinched curvatures in the pointwise sense to the
cases of submanifolds in a Riemannian manifold with codimension
$q(\ge0).$
\medskip

\noindent\textbf{Theorem 3.6.} \emph{Let $M^n$ be an $n$-dimensional
complete submanifold in an $(n + q)$-dimensional Riemannian manifold
$N^{n+q}$. If
$$|A|^2<\frac{8}{3}\Big(\overline{K}_{\min}-\frac{1}{4}\overline{K}_{\max}\Big)+
\frac{|H|^2}{n-1},$$
then $M$ is diffeomorphic to a space form. In particular, if $M$ is simply connected, then $M$ is diffeomorphic
to $\mathbb{S}^n$ or $\mathbb{R}^n$.}\medskip

Let $N^{n+q}$ be an $(n+q)$-dimensional complete Riemannian
manifold. Suppose that $N$ satisfies\\
\indent (1) the sectional curvature $-K_1\leq K_N\leq K_2$ for $K_1,
K_2\geq0$;\\
\indent (2) $|\bar{\nabla}\bar{R}|\leq L$ for $L\geq0$;\\
\indent (3) the injectivity radius $inj(N)\geq i(N)>0$.\\
$N$ is called a Riemannian manifold with bounded geometry.
We call $K_1, K_2, L,$ and $i_N$ the bound constants.
\medskip

In 2012, Liu-Xu-Zhao \cite{liu2012mean} proved the following convergence theorem.

\medskip

\noindent\textbf{Theorem 3.7.} \emph{Let $F_0:M^{n} \rightarrow
N^{n+q}$ be an $n$-dimensional closed submanifold in an
$(n+q)$-dimensional complete Riemannian manifold with bounded
geometry with bound constants $K_1, K_2, L,$ and $i_N$. There is an
explicit constant $b_0$ depending on $n$, $q$, $K_1$, $K_2$ and $L$
such that if $F_0$ satisfies
\begin{eqnarray}
\label{pinch-cond} |A|^2<\begin{cases}
\frac{4}{3n}|H|^2-b_0, \ &n = 2, 3, \\
\frac{1}{n-1}|H|^2-b_0, \ &n \geq 4,
\end{cases}
\end{eqnarray}
then the mean curvature flow with $F_0$
as initial value contracts to a round point in finite time.
In particular, $M$ is diffeomorphic to
$\mathbb{S}^n$.}\medskip

In \cite{LXYZ}, Liu-Xu-Ye-Zhao obtained some convergence theorems
for the mean curvature flow of closed submanifolds in Euclidean
spaces under suitable integral curvature pinching conditions.

\section{MCF with positive Ricci curvature}
There are very few optimal differentiable spheres in curvature and
topology of manifolds. In 1995, Grove and Wilhelm \cite{GW} first
verified an optimal differentiable sphere theorem for a class of
Riemannian manifolds with positive sectional curvature. Using the
Ricci flow and stable currents, Xu and Gu \cite{XG1} proved an
optimal differentiable sphere theorem for certain submanifolds with
positive sectional curvature. A challenging problem is: is it
possible to prove an optimal differentiable sphere theorem for
certain manifolds with positive Ricci curvature?

In 2015, the authors \cite{leiXuOptimal} first proved the following
optimal convergence theorem for the mean curvature flow of arbitrary
codimension in hyperbolic spaces whose initial submanifold has
positive Ricci curvature.

\medskip

\noindent\textbf{Theorem 4.1.}
\emph{Let $F_0:M\rightarrow \mathbb{F}^{n+q}(c)$ be an $n$-dimensional $(n\geq6)$
complete submanifold in a hyperbolic space with negative constant
curvature $c$. If $M$ satisfies
\begin{eqnarray*}
    \sup_{M} (|A|^2 - \alpha(n,|H|,c)) <0 , \,\,\mbox{\ where\ }\,\,
    n^2c+|H|^2>0,\end{eqnarray*} then the mean curvature flow
with  initial value $F_0$
has smooth solution $F_t(\cdot)$, and $F_t(\cdot)$ converges to a
round point in finite time. In particular, $M$ is diffeomorphic to
$\mathbb{S}^n$.}\medskip

\begin{remark}
Since $\alpha ( n, | H | ,c ) > \frac{1}{n-1} | H |^{2} +2c$, our
theorem above improves Theorem 3.5 for $n\ge6$. Note that almost all
initial submanifolds in the convergence results possess positive
curvature. The pinching condition in Theorem 4.1 implies that the
Ricci curvature of the initial submanifold is positive, but does not
imply positivity of the sectional curvature. The following example
shows that the pinching condition in Theorem 4.1 is optimal for
arbitrary $n (\ge 6)$. Therefore, Theorem 4.1 is an optimal
convergence theorem for the mean curvature flow of arbitrary
codimension, which implies the first optimal differentiable sphere
theorem for submanifolds with positive Ricci curvature.
\end{remark}

\textsc{Example.}
Let $\lambda$, $\mu$ be positive constants satisfying $\lambda   \mu =-c$
and $\lambda > \sqrt{-c}$, where $c<0$. For $n \ge 3$, we consider the
submanifold $M=\mathbb{F}^{n-1} ( c+ \lambda^{2} ) \times \mathbb{F}^{1} (
c+ \mu^{2} ) \subset \mathbb{H}^{n+q} ( c )$. Then $M$ is a complete
submanifold with parallel mean curvature, which satisfies $| H | \equiv (
n-1 ) \lambda + \mu >n \sqrt{-c}$ and $|A|^{2} \equiv ( n-1 ) \lambda^{2} +
\mu^{2} = \alpha ( n, | H | ,c )$.\medskip

The key ingredient of the proof of Theorem 4.1 is to establish the
elaborate estimates for the pinching quantity $\mathring{\alpha} =
\alpha (n, | H | ,c ) - \frac{1}{n} |H|^2$, because our pinching
condition is sharper than that in Theorem 3.5. Using the properties
of $\mathring{\alpha}$ and the evolution equations, we first derive
that $|\mathring{A}|^2< \mathring{\alpha}$ is preserved along the
mean curvature flow. Applying a new auxiliary function $f_{\sigma} =
|\mathring{A}|^2/\mathring{\alpha}^{1-\sigma}$, we deduce that
$|\mathring{A}|^2 \le C_0|H|^{2(1-\sigma )}$ via the De Giorgi
iteration. We then obtain an estimate for $|\nabla H |$. Finally,
using estimates for $|\nabla H|$ and the Ricci curvature, we show
that $\mathrm{diam}(M_{t})\rightarrow 0$ and $|H|_{\min}/| H
|_{\max} \rightarrow 1$ as $t\rightarrow T$. This implies the flow
shrinks to a round point.\medskip

\noindent Proof of Theorem 4.1. We split the proof into several steps.

\noindent\textbf{Step 1.}
We first show that the pinching condition in Theorem 4.1 is
preserved under the mean curvature flow.

Suppose that $M_{0}$ is an $n(\ge 6 )$-dimensional  complete submanifold
satisfying
\[
\sup_{M}(| A |^{2} - \alpha ( n, | H | ,c ))<0 ,\qquad | H |^{2}
+n^{2}c>0.
\]
From a theorem due to Shiohama-Xu \cite{Shiohama}, we see that
$M_{0}$ is actually a compact submanifold.

For a positive integer $n$ and a negative constant $c$, we define a function $\mathring{\alpha} : ( -n^{2} c,+ \infty ) \rightarrow   \mathbb{R}$  by
\begin{equation*}
\mathring{\alpha} ( x ) :=n c+ \frac{n^{2} -2n+2}{2 ( n-1 ) n}  x-
\frac{n-2}{2 ( n-1 )} \sqrt{x^{2} +4 ( n-1 ) c x} .
\end{equation*}
Then $ | A |^{2} < \alpha ( n, | H | ,c ) $ is equivalent to $|\mathring{A}|^2<\mathring{\alpha} (  | H |^2 )$.

Note that
$\mathring{\alpha}(|H|^2) \rightarrow 0$  as $
|H |^{2} \rightarrow -n^{2} c.$
This implies that if
$|\mathring{A}|^2<\mathring{\alpha} (  | H |^2 )$ remains true along the flow,
then $| H |^{2} +n^{2} c>0$ also remains true.

We consider $U=|\mathring{A}|^2-\mathring{\alpha}(|H|^2)$. Along the
mean curvature flow, $U$ satisfies
\begin{align*}
&\left(\partial_t-\Delta \right)U\\
=& \left(\partial_t-\Delta \right)|\mathring{A}|^2
-\mathring{\alpha}'(|H|^2)\cdot \left(\partial_t-\Delta \right)|H|^2\\
&+\mathring{\alpha}''(|H|^2) \big|\nabla | H |^{2} \big|^{2}\\
\le& - 2 | \nabla A
|^2 + 2 \left( \tfrac{1}{n} + \mathring{\alpha}' (| H |^2) \right) | \nabla H
|^2 + \mathring{\alpha}'' (| H |^2) \cdot | \nabla | H |^2 |^2 \nonumber\\
&   + 2 | \mathring{A} |^2 \left( | A |^2 - n c \right)
- 2 \mathring{\alpha}' (| H |^2) | H |^2
\left( | A |^2  + n c \right) .
\end{align*}

For convenience, we denote $\mathring{\alpha} ( | H |^{2} )$,
$\mathring{\alpha}' ( | H |^{2} )$ and $\mathring{\alpha}'' ( | H
|^{2} )$ by $\mathring{\alpha}$, $\mathring{\alpha}'$ and
$\mathring{\alpha}''$, respectively.

Using the inequality $| \nabla A |^2 \ge \frac{3}{n + 2} | \nabla H
|^2$ on submanifolds of space forms, we obtain
\begin{align*}
&\left( \partial_t - \Delta \right)U \\
\le & \ 2 U \left( 2 \mathring{\alpha} + \tfrac{1}{n} | H |^2 - n c -
\mathring{\alpha}' \cdot| H |^2 +U \right)\\
&+ 2 \left[ -
\tfrac{2 (n - 1)}{n (n + 2)} + \mathring{\alpha}' + 2 | H |^2
\mathring{\alpha}''  \right] | \nabla H |^2\\
& + 2 \left[ \mathring{\alpha} \cdot \Big( \mathring{\alpha} +
\tfrac{1}{n} | H |^2 - n c \Big) - \mathring{\alpha}' \cdot | H |^2\cdot \Big(
\mathring{\alpha} + \tfrac{1}{n} | H |^2 + n c \Big) \right]  .
\end{align*}
To apply the maximum principle for the parabolic partial differential equation, we need to prove the
two expressions in the square brackets of the above formula are non-positive.

By some computations, we prove the following.
\medskip

\noindent\textbf{Lemma 4.2.}
\emph{If $x>-n^{2} c$, then $\mathring{\alpha}$ has the following properties.}

(i) $2 x  \mathring{\alpha}'' ( x ) + \mathring{\alpha}' ( x ) < \frac{2
    ( n-1 )}{n ( n+2 )},$

(ii) $x  \mathring{\alpha}' ( x ) \cdot \left( \mathring{\alpha} ( x ) +
\tfrac{1}{n} x+n c \right) \equiv \mathring{\alpha} ( x ) \cdot \left(
\mathring{\alpha} ( x ) + \tfrac{1}{n} x-n c \right).$\medskip

It follows from Lemma 4.2 that $\left(\partial_t-\Delta \right)U\leq
0$ at the point where $U=0$. Hence, if $U$ is initially negative,
then $U<0$ is preserved along the mean curvature flow. Therefore, we
obtain the following.
\medskip

\noindent\textbf{Lemma 4.3.}
\emph{If $M_{0}$ satisfies $|\mathring{A}|^{2} < \mathring{\alpha} $
and $| H |^{2} +n^{2} c>0$, then there exist small positive constants
$\varepsilon$ and $\delta$, such that for all time $t \in [ 0,T_{max} )$, we have
\[|\mathring{A}|^{2} < \mathring{\alpha} -
\varepsilon   | H |^{2} \quad and \quad | H |^{2} +n^{2} c>\delta.\]
}

\noindent\textbf{Step 2.}
We next prove a pinching estimate  for the traceless
second fundamental form, which guarantees that $M_t$ becomes
spherical along the mean curvature flow.\medskip

\noindent\textbf{Lemma 4.4.}
\emph{There are positive constants $C_0$ and $\sigma$ independent of $t$
such that
\begin{eqnarray*}
    |\mathring{A}|^2\leq C_0|H|^{2(1-\sigma)}
\end{eqnarray*}
holds along the mean curvature flow.}\medskip

To prove the above lemma, we introduce an auxiliary function:
\[ f_{\sigma} = \frac{|\mathring{A}|^{2}}{\mathring{\alpha}^{1- \sigma}} , \]
where $\sigma\in (0,1)$. We need to prove that $f_{\sigma}$ is bounded.

Computing the evolution equation of $f_{\sigma}$, we get
\[ \frac{\partial}{\partial t} f_{\sigma} \le \Delta f_{\sigma} + \frac{2}{|\mathring{A}|} | \nabla
f_{\sigma} | | \nabla H | -12  \varepsilon \frac{f_{\sigma}}{|\mathring{A}|^{2}}
| \nabla H |^{2} -4 c f_{\sigma} + \sigma | H |^{2} f_{\sigma} . \]

Because of the existence of term $\sigma|H|^2f_\sigma$ in the
inequality above, the maximum principle doesn't work here. To deal
with this term, we need the following estimate.

\medskip

\noindent\textbf{Lemma 4.5.}
    \emph{If a submanifold in $\mathbb{H}^{n+q} ( c )$ satisfies}
    $|\mathring{A}|^{2} < \mathring{\alpha} - \varepsilon  |H|^2$ and $| H |^{2}
    +n^{2} c>0$, then we have
    \[ \Delta |\mathring{A}|^{2} \ge 2 \left\langle \mathring{A} , \nabla^{2} H
    \right\rangle +  \frac{\varepsilon}{2}   | H |^{2}   |\mathring{A}|^{2} . \]

From Lemma 4.5, we obtain the following inequality.

\[
\Delta f_{\sigma}
\ge  \frac{2f_{\sigma}}{| \mathring{A}|^{2}}   \left\langle \mathring{A} ,
\nabla^{2} H \right\rangle
+ \frac{\varepsilon}{2} | H |^{2} f_{\sigma} -  (
1- \sigma ) \frac{f_{\sigma}}{\mathring{\alpha}} \Delta \mathring{\alpha} -
\frac{2}{| \mathring{A}|} | \nabla f_{\sigma} | | \nabla H |  .
\]

Integrating both sides of the inequality above and making use of the
divergence theorem, we get
\[
\frac{\varepsilon}{2} \int_{M_{t}} | H |^{2} f_{\sigma}^{p} \mathrm{d} \mu_{t}
\le \int_{M_{t}} \left( \frac{3p f_{\sigma}^{p-1}}{|\mathring{A}|}   | \nabla
f_{\sigma} | |   \nabla H |  +  \frac{5f_{\sigma}^{p}}{|\mathring{A}|^{2}}   |
\nabla H |^{2} \right)   \mathrm{d} \mu_{t} . \label{H2fsp}
\]

Then we get the following estimate for the time derivative of the
integral of $f_{\sigma}$.
\begin{eqnarray*}
    & &\frac{\mathrm{d}}{\mathrm{d} t} \int_{M_{t}} f_{\sigma}^{p} \mathrm{d} \mu_{t}\\
    & = & p
    \int_{M_{t}} f_{\sigma}^{p-1} \frac{\partial f_{\sigma}}{\partial t} \mathrm{d}
    \mu_{t} - \int_{M_{t}} f_{\sigma}^{p} | H |^{2} \mathrm{d} \mu_{t} \nonumber\\
    & \le & p \int_{M_{t}} f_{\sigma}^{p-2} \left[ - ( p-1 ) | \nabla f_{\sigma}
    |^{2} + \left( 2+ \frac{6  \sigma  p}{\varepsilon} \right)
    \frac{f_{\sigma}}{| \mathring{A}|} | \nabla f_{\sigma} | | \nabla H | \right.
    \nonumber\\
    &  & \left. - \left( 12  \varepsilon - \frac{10  \sigma}{\varepsilon}
    \right) \frac{f_{\sigma}^{2}}{| \mathring{A}|^{2}} | \nabla H |^{2} \right] \mathrm{d}
    \mu_{t} -4 c p \int_{M_{t}} f_{\sigma}^{p} \mathrm{d} \mu_{t} . \nonumber
\end{eqnarray*}

Note that the expression in the square bracket of the above formula is a quadratic  polynomial.
By choosing suitable constants $\sigma$ and $p$, we make this quadratic  polynomial non-positive.

Hence we get
$$ \frac{\mathrm{d}}{\mathrm{d} t} \int_{M_{t}} f_{\sigma}^{p} \mathrm{d} \mu_{t}
\le -4 c p \int_{M_{t}} f_{\sigma}^{p} \mathrm{d} \mu_{t} . $$

It can be proved that the maximal existence time of the mean curvature flow is finite.
Therefore, we obtain

\medskip

\noindent\textbf{Lemma 4.6.}
    \emph{There exist a constant $C$ independent of $t$, such that}
\[
    \bigg(\int_{M_t}f_\sigma^pd\mu_t\bigg)^{\frac{1}{p}}\leq C.
\]

Then we show that $f_{\sigma}$ is bounded along the mean curvature
flow via a Stampacchia iteration procedure.\medskip

\noindent\textbf{Step 3.}
We establish a gradient estimate for the mean curvature flow, which
will be used to compare the mean curvature at different points of
the submanifold.\medskip

\noindent\textbf{Lemma 4.7.} \emph{For every $\eta>0$, there exists
a constant $C_\eta$ independent of $t$ such that for all $t\in
[0,T_{\max})$, the following inequality holds}
\begin{eqnarray*}
    \label{16}|\nabla H|^2\leq \eta |H|^4+C_\eta.
\end{eqnarray*}
\begin{proof}

    Consider the function $$f=|\nabla H|^2+(N_1+N_2|
    A|^2)|\mathring{A}|^2-\eta|H|^4,$$ where $N_1$ and $N_2$ are
    sufficiently large positive constants independent of $t$.

    By the evolution equations of $|H|^4$ and $|\nabla H|^2$, we obtain
    \begin{eqnarray*}
        \frac{\partial}{\partial t}f\leq \triangle f+C_5.
    \end{eqnarray*}

    The assertion follows from the maximum principle and the definition of
    $f$.
\end{proof}

\noindent \textbf{Step 4.}
Now we finish the proof of Theorem 4.1.

    By an  estimate for Ricci curvature due to Shiohama-Xu \cite{Shiohama},
    we have
    \begin{eqnarray*}
        \operatorname{Ric}_{M_t}  \ge  \frac{n-1}{n} \left( n c+ \frac{2}{n} | H
        |^{2} - | A |^{2} - \tfrac{n-2}{\sqrt{n ( n-1 )}} | H | |\mathring{A}| \right).
    \end{eqnarray*}

    By the preserved pinching condition
    $|\mathring{A}|^{2} < \mathring{\alpha} - \varepsilon  |H|^2$
    and the following identity
    $$\tfrac{n-2}{\sqrt{n ( n-1 )}} \sqrt{x  \mathring{\alpha} ( x )}
    \equiv \frac{1}{n} x- \mathring{\alpha} ( x ) +n c,$$
    we see that there exists a positive constant $\varepsilon_0$
    independent of $t$ such that
    $ \operatorname{Ric}_{M_t} \ge \varepsilon_0  | H |^{2} . $

    Since $T_{\max}$ is finite, $\max_{M_t}|A|^2\rightarrow \infty$ as
    $t\rightarrow T_{\max}$.
    At a time $t$ close to $T_{\max}$, let $x$ be a point on $M_{t}$ where $| H |$
    achieves its maximum. From the gradient estimate of $H$, along all geodesics of length $l= ( 2  \eta   | H
    |_{\max} )^{-1}$ starting from $x$, we have $| H | >
     ( 1- \eta ) | H |_{\max}$. With $\eta$ small enough, the Ricci curvature of $M_t$ satisfies $\mathrm{Ric}_{M_t} >
     ( n-1 ) \pi^{2} /l^{2}$ on these geodesics. Then from
    Myers' theorem, these geodesics can reach any point of $M_{t}$.

    Thus we have $| H |_{\min} > ( 1- \eta ) | H |_{\max}$ and $\mathrm{diam}
    (M_{t}) \le ( 2  \eta   | H |_{\max} )^{-1}$. So, $\frac{\max_{M_t}|H|}{\min_{M_t}|H|}\rightarrow 1$
    and $\mathrm{diam}(M_{t})\rightarrow 0$  as
    $t\rightarrow T_{\max}$.
    Hence $M_t$'s converge to a single point $o$ as $t\rightarrow
    T_{\max}$.

    Take a rescaling around $o$
    such that the total area of the expanded submanifolds are fixed.
    Then the rescaled immersions converge to a totally umbilical
    immersion as $t\rightarrow T_{\max}$.
    \medskip

Recently, Lei-Xu \cite{leiXuNew} proved the following sharp
convergence theorem for the mean curvature flow of hypersurfaces in
spheres.\medskip

\noindent\textbf{Theorem 4.8.}
 \emph{Let $F_{0} :M^{n} \rightarrow \mathbb{F}^{n+1} (
c)$ be an n-dimensional ($n \geq 3$) closed hypersurface
immersed in a spherical space form. If $F_{0}$ satisfies
\[ |A|^{2} < \gamma_1 ( n,|H|,c ) , \]
then the mean curvature flow with initial value $F_{0}$ has a unique smooth
solution $F: M \times [ 0,T ) \rightarrow \mathbb{F}^{n+1} (c)$,
and $F_{t}$ converges to a round point in finite time, or
converges to a totally geodesic hypersurface as $t \rightarrow \infty$.
Here $\gamma_1 ( n,H,c )$ is an explicit
positive scalar defined by
\[ \gamma_1 ( n,H,c ) = \min \{ \alpha ( H^{2} ) , \beta ( H^{2} ) \} , \]
where
\[ \alpha ( x ) =n c+ \frac{n}{2 ( n-1 )} x- \frac{n-2}{2 ( n-1 )} \sqrt{x^{2}
    +4 ( n-1 ) c x} , \]
\[ \beta ( x ) = \alpha ( x_{0} ) + \alpha' ( x_{0} ) ( x-x_{0} ) +
\frac{1}{2} \alpha'' ( x_{0} ) ( x-x_{0} )^{2} , \]
\[ x_{0} =y_{n} c, \hspace{1em} y_{n} =4 ( 1-n ) + \frac{2 ( n^{2} -4
    )}{\sqrt{2n-5}}   \cos \left( \frac{1}{3} \arctan \tfrac{n^{2} -4n+6}{2 (
    n-1 ) \sqrt{2n-5}} \right) . \]}\medskip

\begin{remark}
The scalar $\gamma_1 ( n,|H|,c )$ satisfies (i) $\gamma_1 ( n,H,c )
> \frac{1}{n-1} H^{2} +2c$; (ii) $\gamma_1 ( n,|H|,c ) > \frac{9}{5}
\sqrt{n-1} c$; (iii) $\gamma_1 ( n,|H|,c )=\alpha(n,|H|,c)$, for
$|H|^2\ge y_n c$. Thus Theorem 4.8 substantially improves the famous
convergence theorem due to Huisken \cite{Huisken2}.
\end{remark}

Now we state the idea of the proof of Theorem 4.8.\medskip

\noindent\textbf{Step 1.}
We first prove that the pinching condition in Theorem 4.8 is
preserved under the mean curvature flow.

Suppose that $M_{0}$ is an $n$-dimensional $( n \ge 3 )$ closed
hypersurface satisfying $ | A |^{2} < \gamma_1 ( | H |^2 ) .$ Here
$\gamma_1$ is an explicit function defined by
\[
\gamma_1 ( x ) = \left\{\begin{array}{ll}
\alpha ( x ) , & x \geq x_{0} ,\\
\beta ( x ) , & 0 \leq x<x_{0} .
\end{array}\right.
\]

Consider $U=|A|^2-\gamma_1(H^2)$. Computing the evolution equation
of $U$, we get
\begin{eqnarray*}
    &&\left( \frac{\partial}{\partial t} - \Delta \right)U  \\
    & \le & 2 U \left( 2 \gamma_1(H^2)   - n c -
    \gamma_1' (H^2)H^2 +U \right)\\
    &  &+ 2 \left[ -
    \frac{3}{n + 2} + \gamma_1'(H^2) + 2 H^2
    \gamma_1''(H^2)  \right] | \nabla H |^2\\
    &  & +2 \Big[2 c H^{2} +  \gamma_1(H^2)^{2} - n c  \gamma_1(H^2) - ( \gamma_1(H^2) +n c ) H^{2}
    \gamma'_1(H^2)  \Big]
    .
\end{eqnarray*}
To apply the maximum principle, we need to prove the two expressions in the
square brackets of the above formula are non-positive.

By some complicated calculations, we get the following lemma.
\medskip

\noindent\textbf{Lemma 4.9.}
\emph{For $x\ge 0$, the function $\gamma_1$ has the following properties.}

(i)  $2x \gamma_1'' ( x ) + \gamma_1' ( x ) \leq\frac{3}{n+2}$,

(ii) $( \gamma_1 ( x ) +n c ) x  \gamma_1' ( x ) \geq 2 c x+  \gamma_1 ( x
)^{2} -n c  \gamma_1 ( x )$.\medskip

It follows from Lemma 4.9 that $\left(\partial_t-\Delta \right)U\leq
0$ at the point where $U=0$. Applying the maximum principle, we
prove that if $U$ is initially negative, then $U<0$ is preserved
along the mean curvature flow.\medskip

\noindent\textbf{Step 2.}
We prove a pinching estimate  for the traceless
second fundamental form, which guarantees that $M_t$ becomes
spherical along the mean curvature flow.\medskip

\noindent\textbf{Lemma 4.10.}
\emph{There are positive constants $C_0$ and $\sigma$ independent of $t$
    such that
    \begin{eqnarray*}
        |\mathring{A}|^2\leq C_0(H^2+c)^{1-\sigma}\operatorname{e}^{-2\sigma c t}
    \end{eqnarray*}
    holds along the mean curvature flow.}\medskip

To prove the lemma above, we introduce an auxiliary function:
\[ f_{\sigma} = |\mathring{A}|^{2} / \gamma_1(H^2)^{1- \sigma} ,\]
where $\sigma\in (0,1)$. We need to show that $f_{\sigma}$ decays
exponentially.

Let $p$ be a sufficiently large number. Computing the evolution equation of $f_\sigma^p$,
we get
$$ \frac{\mathrm{d}}{\mathrm{d} t} \int_{M_{t}} f_{\sigma}^{p} \mathrm{d} \mu_{t}
\le -3\sigma c p \int_{M_{t}} f_{\sigma}^{p} \mathrm{d} \mu_{t} . $$

Thus we obtain

\medskip

\noindent\textbf{Lemma 4.11.}
\emph{There exist a constant $C$ independent of $t$, such that}
\begin{eqnarray*}\label{8-ineq}
    \bigg(\int_{M_t}f_\sigma^pd\mu_t\bigg)^{\frac{1}{p}}\leq C \operatorname{e}^{-3\sigma c t}.
\end{eqnarray*}

Let $g_\sigma=f_\sigma \operatorname{e}^{2\sigma c t}$. Then we show
that $g_{\sigma}$ is bounded along the mean curvature flow via a
Stampacchia iteration procedure. This proves the lemma.\medskip

\noindent\textbf{Step 3.} We establish a gradient estimate for the
mean curvature flow, which will be used to compare the mean
curvature at different points of the hypersurface.\medskip

\noindent\textbf{Lemma 4.12.}
\emph{For every $\eta>0$, there exists a constant $C_\eta$ independent of
    $t$ such that for all $t\in [0,T_{\max})$, there holds}
\[ | \nabla H |^{2} < [ ( \eta  H )^{4} +C ( \eta )^{2} ] \operatorname{e}^{- \sigma c
    t} . \]

\begin{proof}

    Consider the function $$f=\left[|\nabla H|^2+(N_1+N_2|
    A|^2)|\mathring{A}|^2\right]\operatorname{e}^{ \sigma c
        t}-\eta|H|^4,$$ where $N_1$ and $N_2$ are
    sufficiently large positive constants independent of $t$.
    By the evolution equations of $|H|^4$ and $|\nabla H|^2$, we obtain
    \begin{eqnarray*}
        \frac{\partial}{\partial t}f\leq \triangle f+C_5 \operatorname{e}^{- \sigma c
            t}.
    \end{eqnarray*}
    The assertion follows from the maximum principle and the definition of
    $f$.
\end{proof}

\noindent \textbf{Step 4.}  Now we finish the proof of Theorem 4.8.
We have two cases as follows.\medskip

    Case (i). Suppose that $T_{\max}$ is finite.
    By an estimate for Ricci curvature due to Shiohama-Xu \cite{Shiohama},
    we obtain
    $ \operatorname{Ric}_{M_t} \ge \varepsilon_0  ( H ^{2}+c) . $
    Since $T_{\max}$ is finite, $\max_{M_t}|A|^2\rightarrow \infty$ as
    $t\rightarrow T_{\max}$.
    Using the gradient estimate for $H$ and Myers' theorem, we obtain
    $\frac{\max_{M_t}|H|}{\min_{M_t}|H|}\rightarrow 1$
    and $\mathrm{diam}(M_{t})\rightarrow 0$  as
    $t\rightarrow T_{\max}$.
    Hence $M_t$'s converge to a single point $o$ as $t\rightarrow
    T_{\max}$.

    We take a rescaling around $o$.
    Then the rescaled immersions converge to a totally umbilical
    immersion as $t\rightarrow T_{\max}$.\medskip

    Case (ii). Suppose that $T_{\max}=\infty$.
    It follows from the gradient estimate of $H$ that $ H ^2
    <C \operatorname{e}^{- \sigma c t}$.
    Noting that $|\mathring{A}|^2\leq C_0H^{2(1-\sigma)}\cdot \operatorname{e}^{-2\sigma c t}$, we get $| h |^{2}
    \leq C \operatorname{e}^{- \sigma c t}$.
    Since $| h | \rightarrow 0$ as $t
    \rightarrow \infty$, $M_t$ converges to a totally geodesic hypersurface.
    \medskip

    This completes the proof of Theorem 4.8.\\

More recently, Lei-Xu \cite{leiXuMean} proved a sharp convergence
theorem for mean curvature flow ($n\ge6$) of arbitrary codimension
in spheres, which substantially improves the convergence theorem for
mean curvature flow of arbitrary codimension in spheres due to Baker
\cite{Baker}.\medskip

\noindent\textbf{Theorem 4.13.}
\emph{Let $F_{0} :M^{n} \rightarrow \mathbb{F}^{n+q} (
c)$ be an n-dimensional ($n \geq 6$) closed submanifold
immersed in a spherical space form. If $F_{0}$ satisfies
\[ |A|^{2} < \gamma ( n,|H|,c ) , \]
then the mean curvature flow with initial value $F_{0}$ has a unique smooth
solution $F: M \times [ 0,T ) \rightarrow \mathbb{F}^{n+q} (c)$,
and $F_{t}$ converges to a round point in finite time, or
converges to a totally geodesic sphere as $t \rightarrow \infty$.
Here $\gamma( n,|H|,c )$ is an explicit
positive scalar defined by
$$ \gamma( n,|H|,c )=\min \{ \alpha ( | H |^{2} ) , \beta ( | H
|^{2} ) \}, $$  where
\begin{equation*}
\alpha ( x ) =n c+ \frac{n}{2 ( n-1 )} x- \frac{n-2}{2 ( n-1 )} \sqrt{x^{2}
    +4 ( n-1 ) c x},
\end{equation*}
\begin{equation*}
\beta ( x ) = \alpha ( x_{0} ) + \alpha' ( x_{0} ) ( x-x_{0} ) + \frac{1}{2}
\alpha'' ( x_{0} ) ( x-x_{0} )^{2} , \label{beta}
\end{equation*}
\begin{equation*}  x_{0} =y_n c,\quad y_n= \frac{2  n+2 }{n-4} \sqrt{n-1} \left( \sqrt{n-1} - \tfrac{n-4}{2
    n+2 } \right)^{2}. \end{equation*}}\medskip

\begin{remark}
The scalar $\gamma ( n,|H|,c )$ satisfies the following: (i) $\gamma
( n,H,c )
> \frac{1}{n-1} H^{2} +2c$; (ii) $\gamma ( n,|H|,c ) > \frac{7}{6}
\sqrt{n-1} c$; (iii) $\gamma ( n,|H|,c )=\alpha(n,|H|,c)$, for
$|H|^2\ge z_n c$. Hence Theorem 4.13 substantially improves Theorem
3.4.
\end{remark}

In \cite{PiSi2015}, Pipoli and Sinestrari obtained the following
convergence theorem for mean curvature flow of small codimension in
the complex projective space.\medskip

\noindent\textbf{Theorem 4.14.}
\emph{Let $F_0 : M^n \rightarrow \mathbb{C}\mathbb{P}^{\frac{n +
        q}{2}}$ be a closed submanifold of dimension $n$ and codimension $q$ in the
complex projective space with Fubini-Study metric. Let $F : M^n
\times [0, T) \rightarrow \mathbb{C}\mathbb{P}^{\frac{n + q}{2}}$ be
the mean curvature flow with initial value $F_0$. If $F_0$ satisfies
\[ | h |^2 <B_0(n,q,|H|):= \left\{ \begin{array}{ll}
\frac{1}{n - 1} | H |^2 + 2, &  q = 1 \ \mathrm{and}\  n \geq 5,\\
\frac{1}{n - 1} | H |^2 + \frac{n - 3 - 4 q}{n}, & 2 \leq q < \frac{n - 3}{4},
\end{array} \right. \]
then $F_t$ converges to a round point in finite time, or converges to a
totally geodesic submanifold as $t \rightarrow \infty$. In particular, $M$
is diffeomorphic to either $\mathbb{S}^n$ or $\mathbb{C}\mathbb{P}^{n /
    2}$.}\medskip

Most recently, Lei-Xu \cite{leiXuC} proved a new convergence theorem
for mean curvature flow of arbitrary codimension in complex
projective spaces, which substantially improves the convergence
theorem for mean curvature flow of small codimension in complex
projective spaces due to Pipoli-Sinestrari \cite{PiSi2015}.
Precisely, we prove the following theorem.\medskip

\noindent\textbf{Theorem 4.15.}
\emph{Let $F_0 : M^n \rightarrow \mathbb{C}\mathbb{P}^{\frac{n +
        q}{2}}$ be an n-dimensional closed submanifold in
$\mathbb{C}\mathbb{P}^{\frac{n + q}{2}}$. Let $F : M^n \times [0, T)
\rightarrow \mathbb{C}\mathbb{P}^{\frac{n + q}{2}}$ be the mean
curvature flow with initial value $F_0$. If $F_0$ satisfies
\[ | h |^2 <B(n,q,|H|):= \left\{ \begin{array}{ll}
\varphi (| H |^2), & q = 1 \ \mathrm{and} \ n \geq 3,\\
\frac{1}{n - 1} | H |^2 + 2 - \frac{3}{n}, & 2 \leq q < n - 4,\\
\psi (| H |^2), & q \geq n - 4 \geq 2,
\end{array} \right. \]
then $F_t$ converges to a round point in finite time, or converges to a
totally geodesic submanifold as $t \rightarrow \infty$. In particular, $M$
is diffeomorphic to $\mathbb{S}^n$ or $\mathbb{C}\mathbb{P}^{n /
    2}$.}

\emph{Here $\varphi (| H |^2)$ and $\psi (| H |^2)$ are given by
\[ \varphi (| H |^2) = 2 + a_n + \left( b_n + \tfrac{1}{n - 1} \right)
| H |^2 - \sqrt{b_n^2 | H |^4 + 2 a_n b_n | H |^2}, \]
\[ \psi (| H |^2) = \frac{9}{n^2 - 3 n - 3} + \frac{n^2 - 3 n}{n^3 - 4
    n^2 + 3} | H |^2
 - \tfrac{3 \sqrt{| H |^4 + \frac{2}{n} (n - 1) (n^2 - 3) |
        H |^2 + 9 (n - 1)^2}}{n^3 - 4 n^2 + 3}, \]
where $a_n = 2 \sqrt{(n^2 - 4 n + 3) b_n}, \,\, b_n = \min \left\{
\frac{n - 3}{4
    n - 4}, \frac{2 n - 5}{n^2 + n - 2} \right\}$.}\medskip

\begin{remark}
In particular, (i) if $q=1$ and $n\geq 5$, we have $B( n,q,|H| )
>B_0( n,q,|H| ) $ and $B(
n,q,|H| )
> \sqrt{2 (n - 3)}$; \,\,(ii) if $2 \leq q < \frac{n - 3}{4}$, we have $B( n,q,|H| )-B_0( n,q,|H| )=1+\frac{4 q}{n}$. Therefore, Theorem
4.15 substantially improves Theorem 4.14.
\end{remark}

\section{Problem section}

In this section, we give several unsolved problems on the
convergence theorems for the mean curvature flow and sphere theorems
for submanifolds.

Denote by $\mathbb{F}^{n+q}(c)$ the $(n+q)$-dimensional complete
simply connected space form of constant sectional curvature $c$. Let
$M$ be an $n$-dimensional oriented compact submanifold in
$\mathbb{F}^{n+q}(c)$. Based on the discussions in
\cite{GX,GXXZ,leiXuNew,leiXuMean,LXYZ,LXZsome,MN,Shiohama,Shiohama0,XG5,XHZ},
we present the following conjectures.\medskip

\noindent\textbf{Conjecture 5.1 (Liu-Xu-Ye-Zhao \cite{LXYZ}).}
\emph{Let $F_0:M\rightarrow \mathbb{F}^{n+q}(c)$ be an
$n$-dimensional closed submanifold in $\mathbb{F}^{n+q}(c)$ with
$c>0$.  Suppose $F_0$ satisfies
\begin{eqnarray*}
    |A|^2<\alpha(n,|H|,c),
\end{eqnarray*}
then the mean curvature flow with $F_0$ as initial value has a
unique solution $F:M\times[0,T)\rightarrow \mathbb{F}^{n+q}(c)$ on a
maximal time interval, and
either}\\
\emph{$(i)$ $T<\infty$ and $M_t$ converges to a round point as
$t\rightarrow T$; or}\\
\emph{$(ii)$ $T=\infty$ and $M_t$ converges
to a totally geodesic submanifold in $\mathbb{F}^{n+q}(c)$ as
$t\rightarrow \infty$.}

\emph{In particular, $M$ is diffeomorphic to the standard
$n$-sphere.}\\

\noindent\textbf{Conjecture 5.2 (see also \cite{LXYZ}).} \emph{Let
$F_0:M\rightarrow\mathbb{F}^{n+q}(c)$ be an $n$-dimensional closed
submanifold in $\mathbb{F}^{n+q}(c)$. Then there exists a positive
constant $C_n$ depending only on $n$, such that if $M$ satisfies
\begin{eqnarray*}
    \int_{M}|\mathring{A}|^ndM< C_n,
\end{eqnarray*}
then the mean curvature flow with $F_0$ as initial value has a
unique solution $F:M\times[0,T)\rightarrow \mathbb{F}^{n+q}(c)$ on a
maximal time interval, and
either}\\
\emph{$(i)$ $T<\infty$ and $M_t$ converges to a round
point
as $t\rightarrow T$; or}\\
\emph{$(ii)$ $c>0, \, T=\infty$ and $M_t$
converges to a totally geodesic submanifold in $\mathbb{F}^{n+q}(c)$
as $t\rightarrow \infty$.}

\emph{In particular, $M$ is diffeomorphic to the standard
$n$-sphere.}\\

\noindent\textbf{Conjecture 5.3.} \emph{Let
$F_0:M\rightarrow\mathbb{R}^{n+q}$ be an $n$-dimensional closed
submanifold in $\mathbb{R}^{n+q}$. There exists a positive constant
$C_n$ depending only on $n$, such that if $M$ satisfies
\[
    \int_M|H|^ndM<n^nVol(\mathbb{S}^{n})+C_n,
\]
then the mean curvature flow with $F_0$ as initial value has a
unique solution $F:M\times[0,T)\rightarrow \mathbb{R}^{n+q}$ on a
finite maximal time interval, and $M_t$ converges to a round point
as $t\rightarrow T$. In particular, $M$ is diffeomorphic to the
standard $n$-sphere. }\medskip

A challenging problem is: what is the best pinching constant in
Conjecture 5.3? In particular, we have the following stronger
version of the Willmore theorem verified by Marques and Neves
\cite{MN}.\medskip

\noindent\textbf{Conjecture 5.4.}
\emph{Let $F_0:M^2\rightarrow\mathbb{R}^{3}$ be a
    closed surface in $\mathbb{R}^{3}$. If
    $M$ satisfies
    \[
    \int_M |H|^2 dM < 8 \pi ^2,
    \]
    then the mean
    curvature flow with $F_0$ as initial value has a unique solution
    $F:M\times[0,T)\rightarrow \mathbb{R}^{3}$ on a finite maximal
    time interval, and $M_t$ converges to a round point as $t\rightarrow
    T$.
    In particular, $M$ is diffeomorphic to $\mathbb{S}^2$.
}\\

\noindent\textbf{Conjecture 5.5 (Xu-Gu \cite{XG5}).} \emph{Let
$F_0:M\rightarrow \mathbb{F}^{n+q}(c)$ be an $n(\geq 3)$-dimensional
    closed submanifold in an $(n+q)$-dimensional space form
    $\mathbb{F}^{n+q}(c)$ with $n^2c+|H|^2>0$. If the Ricci curvature of
    $M$ satisfies
    $$Ric_{M} >
    (n-2)(c+\frac{1}{n^2}|H|^2),$$ then the mean curvature flow
    with $F_0$ as initial value has a
    unique solution $F_t(\cdot)$ on a maximal
    time interval, and $M_t$ converges to a
    round point in finite time, or $c>0$ and $M_t$ converges to a
    totally geodesic sphere as $t \rightarrow \infty$. In particular,
    $M$ is diffeomorphic to $\mathbb{S}^n$.}\\

\noindent\textbf{Conjecture 5.6 (Gu-Xu \cite{GX}).} \emph{Let
$F_0:M\rightarrow \mathbb{F}^{n+q}(c)$ be an $n$-dimensional
    closed submanifold in an $(n+q)$-dimensional space form
    $\mathbb{F}^{n+q}(c)$ with $n^2c+|H|^2>0$. Set
    $k(q,n)=\min\{q\cdot sgn(q-1), n\}.$ If the sectional curvature of $M$
    satisfies
    $$K_M >
    \frac{k(q,n)(c+\frac{1}{n^2}|H|^2)}{2[k(q,n)+1]},$$  then the mean
    curvature flow
    with $F_0$ as initial value has a
    unique solution $F_t(\cdot)$ on a maximal
    time interval, and $M_t$ converges to a
    round point in finite time, or $c>0$ and $M_t$ converges to a
    totally geodesic sphere as $t \rightarrow \infty$. In particular,
    $M$ is diffeomorphic to $\mathbb{S}^n$.}\medskip

Even for $q = 1$ and
$c = 1$, the above problem is still open.\medskip

\noindent\textbf{Conjecture 5.7.} \emph{Let $F_0:M\rightarrow
\mathbb{S}^{n+q}$ be an $n$-dimensional
    closed submanifold in an $(n+q)$-dimensional unit sphere.
    If $|h(u, u)|^2 < \frac{1}{3}$ for any $u \in UM$,  then the mean
    curvature flow
    with $F_0$ as initial value has a
    unique solution $F_t(\cdot)$ on a maximal
    time interval, and $M_t$ converges to a
    round point in finite time, or converges to a
    totally geodesic sphere as $t \rightarrow \infty$. In particular,
    $M$ is diffeomorphic to $\mathbb{S}^n$.}\medskip

\noindent\textbf{Conjecture 5.8 (Xu-Huang-Zhao \cite{XHZ}).}
\emph{Let $F_0:M\rightarrow \mathbb{S}^{n+q}$ be an $n$-dimensional
    closed submanifold in an $(n+q)$-dimensional unit sphere. Set
    \[\tau (x)=\max_{u,v \in U_x M,u\perp v}|h(u,u)-h(v,v)|^2. \]
    If $\tau (x)<\frac{4}{3}$ for all $x \in M$,  then the mean
    curvature flow
    with $F_0$ as initial value has a
    unique solution $F_t(\cdot)$ on a maximal
    time interval, and $M_t$ converges to a
    round point in finite time, or converges to a
    totally geodesic sphere as $t \rightarrow \infty$. In particular,
    $M$ is diffeomorphic to $\mathbb{S}^n$.}\medskip

Notice that if $|h(u, u)|^2 < \frac{1}{3}$ for any $u \in U_x M$,
then $\tau (x)<\frac{4}{3}$. This implies that Conjecture 5.8 is
stronger than Conjecture 5.7.

\vspace{1em}

Center of Mathematical Sciences\

Zhejiang University\

Hangzhou 310027\

China\\

E-mail address: lei-li@zju.edu.cn; xuhw@zju.edu.cn


\medskip




\begin{thebibliography}{99}
\bibitem{Alencar}H. Alencar and M. do Carmo, {\it
Hypersurfaces with constant mean curvature in spheres}, Proc. Amer.
Math. Soc., {\bf 120}(1994), 1223-1229.

\bibitem{Andrews-Baker}B. Andrews and C. Baker, {\it Mean curvature flow of pinched submanifolds to
spheres,} J. Differential Geom., {\bf85}(2010), 357-395

\bibitem{Araujo}K. O. Ara\'{u}jo and K. Tenenblat, {\it
On submanifolds with parallel mean curvature vector}, Kodai Math.
J., {\bf 32}(2009), 59-76.

\bibitem{Baker}C. Baker, {\it The mean curvature flow of submanifolds of high codimension}, arXiv:1104.4409v1.

\bibitem{Brendle1}S. Brendle, {\it A general convergence result for the Ricci flow in higher dimensions}, Duke Math. J., {\bf 145}(2008), 585-601.

\bibitem{Brendle3}S. Brendle and R. Schoen, {\it Manifolds with $1/4$-pinched curvature are space
forms,} J. Amer. Math. Soc., {\bf22}(2009), 287-307.

\bibitem{ChenOku}B. Y. Chen and M. Okumura, {\it
Scalar curvature, inequality and submanifold}. Proc. Amer. Math.
Soc., {\bf 38} (1973), 605-608.

\bibitem{ChenLi1} J. Y. Chen and J. Y. Li, {\it Singularity of mean curvature
flow of Lagrangian submanifolds}, Invent. Math., {\bf 156}(2004),
25-51.

\bibitem{Cheng}S. Y. Cheng, P. Li and S. T. Yau, {\it Heat equations on minimal submanifold and their applications}, Amer. J. Math., {\bf106}(1984), 1033-1065.

\bibitem{Chern}S. S. Chern, M. do Carmo and S. Kobayashi, {\it Minimal submanifolds of a sphere with second fundamental form of
    constant length,} in Functional Analysis and Related Fields,
Springer-Verlag, New York(1970).

\bibitem{Ejiri}N. Ejiri, {\it Compact minimal submanifolds of a sphere with positive
Ricci curvature,} J. Math. Soc. Japan, {\bf31} (1979), 251-256.

\bibitem{Ge} J. Q. Ge and Z. Z. Tang, {\it A proof of the DDVV conjecture and its equality case},  Pacific J. Math.,
 {\bf237}(2008), 87-95.

\bibitem {GW} K. Grove and F. Wilhelm, {\it Hard and soft packing radius theorems}, Ann. of Math., {\bf 142}(1995), 213-237.

\bibitem{GX} J. R. Gu and H. W. Xu, {\it On Yau rigidity theorem for minimal
submanifolds in spheres}, Math. Res. Lett., {\bf19}(2012),
511-523.

\bibitem{GX1} J. R. Gu and H. W. Xu, {\it The sphere theorems for
manifolds with positive scalar curvature}, J. Differential Geom.,
{\bf92}(2012), 507-545.

\bibitem{XG6} J. R. Gu and H. W. Xu, {\it A sharp differentiable pinching theorem for submanifolds in space forms}, Proc. Amer. Math. Soc., {\bf144}(2016), 337-346.

\bibitem{GXXZ} J. R. Gu, H. W. Xu, Z. Y. Xu and E. T. Zhao, {\it
A survey on rigidity problems in geometry and topology of
submanifolds}, Proceedings of the 6th International Congress of
Chinese Mathematicians,  ALM 37, pp. 79-99, 2016.

\bibitem {Hamilton} R. Hamilton, {\it Three manifolds with positive Ricci curvature}, J. Differential Geom., {\bf 17}(1982), 255-306.

\bibitem{MR772132}G. Huisken, {\it Flow by mean curvature
    of convex surfaces into spheres}, J. Differential Geom., \textbf{20}(1984),
237-266.

\bibitem{MR837523}G. Huisken, {\it Contracting convex
    hypersurfaces in Riemannian manifolds by their mean curvature}, Invent.
Math., \textbf{84}(1986), 463-480.

\bibitem{Huisken2} G. Huisken, {\it Deforming hypersurfaces
of the sphere by their mean curvature,} Math. Z., {\bf195} (1987),
205-219.

\bibitem{Itoh2}T. Itoh, {\it On Veronese manifolds,}  J. Math. Soc. Japan., {\bf 27}(1975), 497-506.

\bibitem{Lawson}B. Lawson, {\it Local rigidity theorems for minimal hypersurfaces,} Ann. of Math., {\bf 89}(1969), 187-197.

\bibitem{Lawson2}B. Lawson and J. Simons, {\it On stable currents and their
    application to global problems in real and complex geometry, } Ann.
of Math., {\bf98}(1973), 427-450.

\bibitem{leiXuOptimal} L. Lei and
H. W. Xu, {\it An optimal convergence theorem for mean
    curvature flow of arbitrary codimension in hyperbolic spaces}, arXiv:1503.06747, 2015.

 \bibitem{leiXuNew}L. Lei and
 H. W. Xu, {\it A new version of Huisken's convergence theorem for mean
    curvature flow in spheres}, arXiv:1505.07217, 2015.

\bibitem{leiXuMean}L. Lei and
H. W. Xu, {\it Mean curvature flow of arbitrary codimension in
    spheres and sharp differentiable sphere theorem}, arXiv:1506.06371v2, 2015.

\bibitem{leiXuC}L. Lei and
H. W. Xu, {\it Mean curvature flow of arbitrary codimension in
complex projective spaces}, arXiv:1605.07963, 2016.

\bibitem{Li0} A. M. Li and J. M. Li, {\it An inequality for matrices and its applications in differential geometry,} Adv. Math.(China), {\bf 20}(1991), 375-376.

\bibitem{Li} A. M. Li and J. M. Li, {\it An intrinsic rigidity theorem for minimal submanifolds in a sphere,} Arch. Math., {\bf 58}(1992), 582-594.

\bibitem{Li2} H. Z. Li, {\it Curvature pinching for odd-dimensional minimal submanifolds in a
sphere},  Publ. Inst. Math. $($Beograd$)$, {\bf53}(1993), 122-132.

\bibitem{Lin}J. M. Lin and C. Y. Xia, {\it Global pinching theorem for even
dimensional minimal submanifolds in a unit sphere}, Math. Z.,
{\bf201}(1989), 381-389.

\bibitem{LXYZ} K. F. Liu, H. W. Xu, F. Ye and E. T. Zhao, {\it The extension and convergence of mean curvature flow in higher
codimension,} Trans. Amer. Math. Soc., {\bf 370}(2018), 2231-2262.

\bibitem{MR3078951}K. F. Liu, H. W. Xu, F. Ye and E. T. Zhao, {\it Mean
    curvature flow of higher codimension in hyperbolic spaces}, Comm. Anal.
Geom., \textbf{21}(2013),  651-669.

\bibitem{liu2012mean}K. F. Liu,
H. W. Xu and E. T. Zhao, {\it Mean
    curvature flow of higher codimension in Riemannian manifolds}, arXiv:1204.0107, 2012.

\bibitem{LXZsome} K. F. Liu, H. W. Xu and E. T. Zhao, {\it Some recent results on mean curvature flow of arbitrary
codimension}, Proceedings of the 6th International Congress of
Chinese Mathematicians, ALM 37, pp. 115-131, 2016.

\bibitem{Lu4} Z. Lu, {\it Normal scalar curvature conjecture and its applications},  J. Funct. Anal., {\bf 261}(2011), 1284-1308.

\bibitem{MN}F. C. Marques and A. Neves, {\it Min-max theory and the Willmore conjecture}, Ann. of Math., \textbf{179}(2014), 683-782.

\bibitem{Okumura}M. Okumura, {\it Submanifolds and a pinching problem on the second fundamental tensor}, Trans. Amer. Math. Soc., {\bf178}(1973), 285-291.

\bibitem{Okumura1}M. Okumura, {\it Hypersurfaces and a pinching problem on the
second fundamental tensor}, Amer. J. Math., {\bf96}(1974),
207-213.

\bibitem{PiSi2015} G. Pipoli and C. Sinestrari, {\it Mean curvature flow of pinched submanifolds of $\mathbb{CP}^n$}, arXiv:1502.00519v3, 2015.

\bibitem{MR1289187} W. Santos, {\it Submanifolds with parallel mean curvature vector in spheres}, Tohoku Math. J., {\bf46}(1994),
   403-415.

\bibitem{Shen}C. L. Shen, {\it A global pinching theorem for minimal
hypersurfaces in sphere}, Proc. Amer. Math. Soc.,
{\bf105}(1989), 192-198.

\bibitem{Shen2}Y. B. Shen, {\it Submanifolds with nonnegative sectional curvature}, Chinese Ann. Math. SerB, {\bf 5}(1984), 625-632.

\bibitem{ShYB} Y. B. Shen, {\it Curvature pinching for three-dimensional minimal
submanifolds in a sphere}, Proc. Amer. Math. Soc.,
{\bf115}(1992), 791-795.

\bibitem{Shiohama}K. Shiohama and H. W. Xu, {\it The topological sphere theorem for complete submanifolds,
} Compositio Math., {\bf107}(1997), 221-232.

\bibitem{Shiohama0} K. Shiohama and H. W. Xu, {\it Rigidity and sphere theorems for
submanifolds}, Kyushu J. Math. I, {\bf 48}(1994), 291-306; II, {\bf
54}(2000), 103-109.

\bibitem{Simon}U. Simon, {\it
Submanifolds with parallel mean curvature vector and the curvature
of minimal submanifolds of spheres}, Arch. Math., (Basel) {\bf
29}(1977), 106-112.

\bibitem{Simons}J. Simons, {\it Minimal varieties in Riemannian manifolds,} Ann. of Math., {\bf 88}(1968), 62-105.

\bibitem{Smoczyk4} K. Smoczyk, {\it Mean curvature flow in higher codimension: introduction and survey}, Global differential geometry, 231-274, Springer Proc. Math., 17, Springer, Heidelberg, 2012.

\bibitem{Wang1} M. T. Wang, {\it Mean curvature flow of surfaces in einstein four-manifolds}, J. Differential Geom.,
{\bf 57}(2001), 301-338.

\bibitem{Wang2} M. T. Wang, {\it Long time existence and convergence of graphic mean curvature flow in arbitrary
codimension}, Invent. Math., {\bf 148}(2002), 525-543.

\bibitem{Wang3} M. T. Wang, {\it Gauss maps of the mean curvature flow}, Math. Res. Lett., {\bf 10}(2003), 287-299.

\bibitem{Wang4} M. T. Wang, {\it A convergence result of the lagrangian mean curvature flow}, Third International Congress of Chinese Mathematicians. Part 1, 2, 291-295,
AMS/IP Stud. Adv. Math., 42, pt. 1, 2, Amer. Math. Soc., Providence,
RI, 2008.

\bibitem{Wang5} M. T. Wang, {\it Lectures on mean curvature flows in higher codimensions}, Handbook of geometric analysis. No. 1, 525-543,
Adv. Lect. Math. (ALM), 7, Int. Press, Somerville, MA, 2008.

\bibitem{Xin}Y. L. Xin, {\it Mean curvature flow with bounded Gauss image}, Results Math., {\bf 59}(2011), 415-436.

\bibitem{Xu90}H. W. Xu, {\it Pinching theorems, global pinching theorems and
eigenvalues for Riemannian
  submanifolds,} Ph.D. dissertation, Fudan University, 1990.

\bibitem{Xu}H. W. Xu, {\it A rigidity theorem for submanifolds with parallel mean curvature in
a sphere,} Arch. Math., {\bf61}(1993), 489-496.

\bibitem {xuhwrigidity}H. W. Xu, {\it Rigidity of submanifolds with parallel mean curvature in space forms}, preprint, 1993.

\bibitem{XG1} H. W. Xu and J. R. Gu, {\it An optimal differentiable sphere theorem for complete
manifolds,} Math. Res. Lett., {\bf17}(2010), 1111-1124.

\bibitem{XG5} H. W. Xu and J. R. Gu, {\it Geometric, topological and differentiable rigidity of
submanifolds in space forms},  Geom. Funct. Anal., {\bf 23}(2013),
1684-1703.

\bibitem{XHZ} H. W. Xu, F. Huang and  E. T. Zhao, {\it
Geometric and differentiable rigidity of submanifolds in spheres},
J. Math. Pures Appl., {\bf 99}(2013), 330-342.

\bibitem{XT} H. W. Xu and L. Tian, {\it A differentiable sphere theorem
inspired by rigidity of minimal submanifolds}, Pacific J. Math., {\bf 254}(2011), 499-510.

\bibitem{XZ}H. W. Xu and E. T. Zhao, {\it Topological and
differentiable sphere theorems for complete submanifolds}, Comm.
Anal. Geom., {\bf 17}(2009), 565-585.

\bibitem{Yau}S. T. Yau, {\it Submanifolds with constant mean curvature I, II,} Amer. J. Math., {\bf96, 97}(1974, 1975), 346-366, 76-100.

\end{thebibliography}
\end{document}